\documentclass{amsart}

\theoremstyle{definition}

\theoremstyle{remark}

\begin{document}

\title{Rational Landen transformations on $\mathbb{R}$}

\author{Dante Manna}
\address{Department of Mathematics and Statistics,
Dalhousie University, Nova Scotia, Canada B3H 3J5}
\email{dmanna@mathstat.dal.ca}

\author{Victor H. Moll}
\address{Department of Mathematics,
Tulane University, New Orleans, LA 70118}
\email{vhm@math.tulane.edu}

\subjclass{Primary 33}

\date{\today}


\begin{abstract}
The Landen transformation $(a,b) \mapsto ( (a+b)/2,
\sqrt{ab} )$ preserves the value of an elliptic integral and its iteration
produces the classical arithmetic-geometric mean $\text{AGM}(a,b)$. We
present analogous transformations for rational functions integrated over
the whole real line.
\end{abstract}

\maketitle

\newcommand{\nn}{\nonumber}
\newcommand{\ba}{\begin{eqnarray}}
\newcommand{\ea}{\end{eqnarray}}
\newcommand{\E}{{\mathfrak{E}}}
\newcommand{\F}{{\mathfrak{F}}}
\newcommand{\Ro}{{\mathfrak{R}}}
\newcommand{\ift}{\int_{0}^{\infty}}
\newcommand{\ifft}{\int_{- \infty}^{\infty}}
\newcommand{\no}{\noindent}
\newcommand{\X}{{\mathbb{X}}}
\newcommand{\Q}{{\mathbb{Q}}}
\newcommand{\R}{{\mathbb{R}}}
\newcommand{\Y}{{\mathbb{Y}}}
\newcommand{\Ftwo}{{{_{2}F_{1}}}}
\newcommand{\realpart}{\mathop{\rm Re}\nolimits}
\newcommand{\imagpart}{\mathop{\rm Im}\nolimits}

\newtheorem{Definition}{\bf Definition}[section]
\newtheorem{Thm}[Definition]{\bf Theorem}
\newtheorem{Example}[Definition]{\bf Example}
\newtheorem{Lem}[Definition]{\bf Lemma}
\newtheorem{Note}[Definition]{\bf Note}
\newtheorem{Cor}[Definition]{\bf Corollary}
\newtheorem{Prop}[Definition]{\bf Proposition}
\newtheorem{Problem}[Definition]{\bf Problem}
\numberwithin{equation}{section}

\section{Introduction} \label{sec-intro}
\setcounter{equation}{0}

The problem of indefinite integration of rational
functions $R(x) = B(x)/A(x)$ was
finished by J. Bernoulli in the eighteenth century. He completed the original
attempt by Leibniz of a general partial decomposition of $R(x)$.
The result is that a primitive of a rational function
is always elementary: it consists of a new rational function (its {\em rational
part}) and the logarithm of a second rational function (its {\em transcendental
part}).

In the middle of the nineteenth century Hermite \cite{hermite1} and
Ostrogradsky \cite{ostro1} developed algorithms to compute the rational part
of the primitive of $R(x)$ {\em without} factoring $A(x)$. More recently
Horowitz \cite{horowitz1} rediscovered this method and discussed its
complexity. The problem of computing the transcendental part of the primitive
was finally solved by Lazard and Rioboo \cite{lazard1}, Rothstein
\cite{rothstein3} and Trager \cite{trager2}. For detailed descriptions and
proofs of these algorithms the reader is referred to \cite{bronstein2} and
\cite{geddes3}.

This paper contains a method of computing {\em definite} rational integrals
that, unlike the methods described above, does not involve the
factorization of any polynomial. In this new method, the value of the integral
is obtained as the limit of a sequence of transformations of the
coefficients of the integrand. Thus, the
algorithm presented here is in the spirit of the classical
{\em Landen transformation} for elliptic integrals. These are integrals of
the form
\ba
K(k) & = &
\int_{0}^{1} \frac{dx}{\sqrt{(1-x^{2})(1-k^{2}x^{2})}},
\label{intK}
\ea
\no
that have been studied since the eighteenth century. The reader will find
in \cite{mckmoll} more information about them. Its trigonometric version, 
\ba
G(a,b) & = & \int_{0}^{\pi/2} \frac{d \theta}{\sqrt{a^{2} \cos^{2} \theta
+ b^{2} \sin^{2} \theta}}, \label{intG}
\nn
\ea
\no
was considered by Gauss \cite{gauss1} in his work on the lemniscate. The
special case $k = i$,
\ba
\text{Lem} & = & \int_{0}^{1} \frac{dx}{\sqrt{1-x^{4}}},
\label{lemniscate}
\ea
\no
appears as an expression for its arc length. He
inferred from a numerical evaluation that
the function $G(a,b)$ is invariant under
\ba
\mathfrak{E}: \, (a,b) & \mapsto \left( \tfrac{1}{2}(a+b), \sqrt{ab} \right).
\label{mapE}
\ea
\no
A transformation of the parameters of an integral is called a {\em Landen
transformation} if it preserves the value of the integral. The 
example (\ref{mapE}) is the original one.

It is a classical result that the iteration of $\mathfrak{E}$ produces two
sequences, $a_{n}$ and $b_{n}$,
that converge quadratically to a common limit: $\text{AGM}(a,b)$, the
{\em arithmetic-geometric mean} of $a$
and $b$. The invariance of the elliptic integral (\ref{intG}) yields
\ba
G(a,b) & = & \frac{\pi}{2 \, \text{AGM}(a,b)}.
\ea
\no
Iteration of (\ref{mapE}) provides a method to evaluate the elliptic
integral $G(a,b)$.
For instance, four steps starting at $a_{0} =1, \,
b_{0} = \sqrt{2}$ yield $22$ correct digits of the  integral
in (\ref{lemniscate}). See \cite{borwein1} for details and its relation to
modern evaluations of $\pi$.

We consider here the space of rational functions
\ba
{\mathfrak{R}}_{p} & := &
\left\{ R(x) = \frac{B(x)}{A(x)} {\Bigg{|}} \,
A(x) = \sum_{k=0}^{p} a_{k}x^{p-k} \text{ and }
B(x) = \sum_{k=0}^{p-2} b_{k}x^{p-2-k}  \, \right\}.
\nn
\ea
\no
We assume
\begin{itemize}
\item{The degree $p$ is an even positive integer.}
\item{The coeffients $a_{k}$ and $b_{k}$ are real numbers.}
\item{The polynomial $A(x)$ has no real zeros.}
\end{itemize}
\no
Under these assumptions the integral
\ba
I & := & \ifft R(x) \, dx
\ea
\no
is finite.

We describe a transformation on  the parameters
\ba
{\mathfrak{P}}_{p} & := & \{ a_{0}, \, a_{1}, \cdots, a_{p}; \, b_{0}, \,
b_{1}, \cdots, b_{p-2} \, \}
\label{initialpara}
\ea
\no
of $R \in {\mathfrak{R}}_{p}$ that preserves the integral $I$. In fact, we
produce a family of maps, indexed by $m \in \mathbb{N}$,
$$
{\mathfrak{L}}_{m,p}: {\mathfrak{R}}_{p} \to {\mathfrak{R}}_{p},
$$
\no
such that
\ba
\ifft R(x) \, dx & = & \ifft {\mathfrak{L}}_{m,p}(R(x)) \, dx.
\ea
\no
The maps ${\mathfrak{L}}_{m,p}$ induce a {\em rational Landen
transformation} on the coefficients:
\ba
\Phi_{m,p}: {\mathbb{R}}^{2p} \to {\mathbb{R}}^{2p}.
\ea
\no
We provide numerical
evidence that the iterates of this map converge to a limit, with convergence
of order $m$. 

In the case $m=p=2$, we will show that the integral
\ba
I(a_{0},a_{1},a_{2}) & = & \ifft \frac{dx}{a_{0}x^{2}+a_{1}x + a_{2}}
\label{example0}
\ea
\no
is invariant under the transformation
\ba
a_{0} & \mapsto & \frac{2a_{0}a_{2}}{a_{0}+a_{2}}, \label{map0} \\
a_{1} & \mapsto & \frac{a_{1}(a_{2}-a_{0})}{a_{0}+a_{2}}, \nn \\
a_{2} & \mapsto & \frac{(a_{0}+a_{2})^{2} - a_{1}^{2}}{2(a_{0}+a_{2})}. \nn
\ea
\no
This example is
discussed in detail in \cite{manna-moll1}.

\medskip

The theory of Landen transformations for rational integrands is
divided into two cases, according to the domain of integration. \\

\no
{\bf Case 1}: The interval of integration is not the whole real line.  \\

Integration over a  finite interval $[a,b]$ is transformed to the half-line
$[0, \infty)$ by a bilinear transformation. In detail,
\ba
\int_{a}^{b} R(x) \, dx & = &
(b-a) \ift R \left( \frac{a + bt}{1+t} \right) \frac{dt}{(1+t)^{2}}.
\ea
\no
Similarly, integration over half-lines $[a, \infty)$ and $(-\infty, a]$ can
be reduced to $[0, \infty)$ by translations and reflections. Thus, the
interval $[0, \infty)$ encompasses all integrals that fall in this case.

Landen transformations for {\em even} rational functions on $[0, \infty)$ were
established in \cite{boros1}. 
For example, the integral
\ba
U_{6}(a,b;c,d,e) & := & \ift \frac{cx^{4} + dx^{2} + e}{x^{6}+ax^{4}+bx^{2}
+1 } \, dx
\label{u6}
\ea
\no
is invariant under
\ba
a & \mapsto & \frac{ab + 5a + 5b + 9}{(a+b+2)^{4/3}} \label{landen6} \\
b & \mapsto & \frac{a + b + 6}{(a+b+2)^{2/3}}, \nn
\ea
\no
with similar rules  for the coefficients $c, \, d$ and $e$.

The map (\ref{landen6}) can be iterated to produce a sequence
$(a_{n},b_{n};c_{n},d_{n},e_{n})$ with the property
\ba
U_{6}(a_{n},b_{n};c_{n},d_{n},e_{n}) & = & U_{6}(a,b;c,d,e).
\ea
\no
Its convergence was discussed in \cite{boros2}, assuming that
the initial conditions $a_{0}, \, b_{0}$ are nonnegative. The 
main result is the
existence of a number $L$, depending on the initial data 
$a_{0}, \cdots, e_{0}$, such that
$a_{n} \to 3, \, b_{n} \to 3,$ and $c_{n} \to L, \, d_{n} \to 2L,$ and $e_{n}
\to L$. The convergence is quadratic.

The positivity  condition on initial data was eliminated
in \cite{hubbard1}, where
we reinterpret the Landen transformation (\ref{landen6}) in geometric terms.
The new integrand is the direct image of the original one under the
map $w = (z^{2}-1)/2z$. In concrete terms, if $R$ is the original
integrand and
\ba
z_{\pm}(w) & = & w \pm \sqrt{w^{2} + 1}
\ea
\no
are the two branches of the inverse of $w$, then the new integrand is given by
\ba
R(z_{+}(w)) \frac{dz_{+}}{dw} +
R(z_{-}(w)) \frac{dz_{-}}{dw}.
\label{newinteg}
\ea
\no
This geometric interpretation extends to the
algorithm presented in \cite{boros1}, where an analogue of (\ref{landen6})
is given for an arbitrary even function.
These transformations on
the coefficients define a map,
\begin{equation}
\Phi_{2n}: {\mathbb{R}}^{2n-1} \to {\mathbb{R}}^{2n-1},
\end{equation}
\no
which is the rational analogue of $\mathfrak{E}$ in (\ref{mapE}). These are
named {\em even rational Landen}.
Using this approach, we have established a proof that the iterations of
$\Phi_{2n}$ converge precisely when the initial
integral is finite.

A purely dynamical proof of convergence of the iterations of $\Phi_{2n}$
is presented in \cite{marc-moll}, but
only for the case of degree $6$. The relation
between (\ref{landen6}) and the
invariance of the rational integrals is still part of the argument. It is
established that the iterations are eventually mapped to the 
first quadrant, and then the results of \cite{boros2} are applied. It
would be desirable to obtain a proof of convergence
completely independent of the integrals that
gave origin to these maps.

\medskip

The existence of this type of transformation for an odd rational integrand
is an open question. 

\medskip

\no
{\bf Case 2}: the domain of integration is the real line. \\

This is the case we present here. We give a
Landen transformation for integrals over $\mathbb{R}$. The
convergence of the iterations of these maps can be established along the
lines of \cite{hubbard1}, but a more direct analysis is still an open question.
The issue of convergence is not discussed here, except for the numerical
examples in Section \ref{sec-landenex}.

The new integrands, $\mathfrak{L}_{m,p}(R(x))$, depend on the parameter $p$, the
degree of the denominator of the original integrand, and the parameter $m$, the
order of convergence. {\em Both parameters are arbitrary}.

Section 2 presents a preliminary
example that illustrates the methods developed in the
rest of the paper. Section 3 introduces two families of polynomials that are
the basis of the rational Landen transformations. Section 4
consists of some simple trigonometrical identities.  The 
integrand is scaled in Section
5, using the polynomials studied in Section 3. The algorithm leading to the
rational Landen transformation is a consequence of the vanishing of
a class of integrals. This is presented in Section 6. Examples are given in
the last section.

\medskip

\section{An example} \label{sec-example}
\setcounter{equation}{0}

We begin with an
example of a Landen transformation
that introduces the methods described in later sections.

The integral of the rational function
\ba
R(x) & = & \frac{x^{2}+x+1}{x^{4}+6x^{3}+29x^{2}+60x+100}
\ea
\no
is evaluated as
\ba
I := \ifft R(x) \, dx & = & \frac{38 \, \pi}{31 \, \sqrt{31}},
\ea
\no
using the factorization
\begin{equation}
x^{4}+6x^{3}+29x^{2}+60x+100 = (x^{2}+3x+10)^{2}.
\end{equation}

We will produce a new rational function,
\ba
{\mathfrak{L}}_{2,4}(R(x)) & = &
\frac{202x^{2}+45x+97}
{400x^{4}+1080x^{3}+2969x^{2}+3024x+3136},
\ea
\no
and show that it satisfies
\ba
\ifft {\mathfrak{L}}_{2,4}(R(x)) \, dx & = & \ifft R(x) \, dx.
\ea
\no
(The notation ${\mathfrak{L}}_{2,4}$ indicates the degrees of the transformation
used to produce this new function. Details are given in Section
\ref{sec-reduction}).

\medskip

The first step is to multiply the denominator,
\ba
A(x) & = & x^{4} + 6x^{3} + 29x^{2} + 60x + 100,
\ea
\no
by
\ba
Z(x) & = & 1600x^{4} - 960x^{3} + 464x^{2} - 96 x  + 16,
\ea
\no
so that $E(x) = A(x) Z(x)$ can be written as a
homogeneous polynomial in the variables
\ba
P_{2}(x) = x^{2}-1 & \text{ and } & Q_{2}(x) = 2x.
\ea
\no
(These polynomials will be described in Section \ref{sec-poly1}.) In detail,
\ba
E(x) & = &  \sum_{l=0}^{4} e_{l} P_{2}^{4-l}(x) Q_{2}^{l}(x),
\ea
\no
with $
e_{0}=1600, \, e_{1}= 4320, \, e_{2} = 11876, \,
e_{3}=12096, \text{ and } e_{4}=12544. $ Then, with
$C(x) = B(x)Z(x)$, we obtain
\ba
I & = & \ifft \frac{B(x)}{A(x)} \, dx  =
\ifft \frac{C(x)}{E(x)} \, dx.
\label{int51}
\ea

Now write
\ba
E(x) & = & Q_{2}^{4}(x) \left( \sum_{l=0}^{4} e_{l} R_{2}(x)^{4-l} \right),
\ea
\no
where
\ba
R_{2}(x) & = & \frac{P_{2}(x)}{Q_{2}(x)} = \frac{x^{2}-1}{2x}.
\ea
\no
We would like to make the change of variables
$y = R_{2}(x)$ in (\ref{int51}). The function $R_{2}(x)$ has
a multivalued inverse,
with its two branches given by
\ba
x & = & y \pm \sqrt{y^{2}+1}.
\ea
\no
Therefore, we must split the evaluation of the original integral at the
singularity $x=0$
of $R_{2}(x)$. The identity (\ref{int51}) is written as
\ba
I & = & \int_{-\infty}^{0} \frac{C(x)}{E(x)} \, dx +
\int_{0}^{\infty} \frac{C(x)}{E(x)} \, dx  \nn \\
 & = & \ifft \frac{N_{-}(y)}{E_{1}(y)} \, dy +
 \ifft \frac{N_{+}(y)}{E_{1}(y)} \, dy  \nn
\ea
\no
where
\ba
E_{1}(y) & = & \sum_{l=0}^{4} e_{l}y^{4-l}
\ea
\no
and
\ba
N_{\pm}(y) & = &
\frac{C(y \pm \sqrt{y^{2}+1})}{Q_{2}^{4}(y \pm \sqrt{y^{2}+1})} \,
\, \frac{d}{dy} \left( y \pm \sqrt{y^{2}+1} \right).
\ea
\no
The new integrand, $(N_{+}(y) + N_{-}(y))/E_{1}(y)$, corresponds to
the expression in (\ref{newinteg}). A direct calculation shows that
\ba
N_{-}(y) + N_{+}(y) & = & 4(202y^{2} + 45y + 97),
\ea
\no
so that
\ba
I & = & \ifft \frac{202y^{2}+45y+97}
{400y^{4}+1080y^{3}+2969y^{2}+3024y+3136} \, dy \\
 & = & \ifft \frac{202y^{2}+45y+97}{(20y^{2}+27y+56)^{2}} \, dy, \nn
\ea
\no
as claimed.

A proof of a transformation of this type for a general rational integrand is
provided in the next four sections.

\section{A family of polynomials} \label{sec-poly1}
\setcounter{equation}{0}

For $m \in \mathbb{N}$, we introduce the polynomials
\ba
P_{m}(x) & = & \sum_{j=0}^{\lfloor{m/2 \rfloor}} (-1)^{j}
\binom{m}{2j} x^{m-2j} \quad \text{ and } \label{polyP} \\
Q_{m}(x) & = & \sum_{j=0}^{\lfloor{(m-1)/2 \rfloor}} (-1)^{j}
\binom{m}{2j+1} x^{m-(2j+1)}, \label{polyQ}
\ea
\no
which play a fundamental role in the algorithm discussed here. They
will comprise the numerators and denominators of a natural change of
variables discussed in the last two sections.

The degrees of $P_{m}$ and
$Q_{m}$ are $m$ and $m-1$, respectively.  Observe that
\ba
P_{2}(x) = x^{2}-1 & \text{ and } & Q_{2}(x)  = 2x
\ea
\no
have appeared in Section \ref{sec-example}.

\medskip

\begin{Prop}
\label{rofm}
Let $M(x) = \frac{x+i}{x-i}$ and $f_{m}(x) = x^{m}$. Then the rational
function $R_{m} = P_{m}/Q_{m}$ satisfies
\ba
R_{m} & = & M^{-1} \circ f_{m} \circ M. \label{relation1}
\ea
\end{Prop}
\begin{proof}
The identity follows from
\ba
(x + i)^{m} + (x-i)^{m}  =  2 P_{m}(x) & \text{ and } &
(x + i)^{m} - (x-i)^{m}  =  2i Q_{m}(x). \nn
\ea
\end{proof}

\medskip

\begin{Cor}
The function $R_{m}$ satisfies
\ba
R_{m}(\cot \theta ) & = & \cot (m \theta).
\label{cotangent}
\ea
\end{Cor}
\begin{proof}
Use $M( \cot \theta) = e^{2 i \theta}$ in (\ref{relation1}).
\end{proof}

\medskip

\begin{Note}
The multiplicative property $R_{n} \circ R_{m} = R_{nm}$ shows
that the functions $R_{m}$ form a family of {\em commuting}
rational functions. The cotangent function in (\ref{cotangent}) appears as 
the limiting case of the Weierstrass elliptic
$\mathfrak{p}$-function,
 \ba \mathfrak{p}(x) & = & \frac{1}{x^{2}}+
\sum_{n_{1}, \, n_{2} \in \mathbb{Z}} \left[ \frac{1}{(x -
n_{1}\omega_{1} -n_{2} \omega_{2})^{2}} -
\frac{1}{(n_{1}\omega_{1} +n_{2} \omega_{2})^{2}}  \right], \ea
\no where the term $n_{1} = n_{2} =0$ is excluded from the sum. In
the case $ \omega_{1} =1$ and $\omega_{2} \to \infty$, we get \ba
\mathfrak{p}(x) \to -\pi \frac{d}{dx} \cot( \pi x) -
\frac{\pi^{2}}{3}. \ea \no The function $\mathfrak{p}(nx)$ is even
and elliptic, therefore it is a rational function $g_n$ of
$\mathfrak{p}$. In view of $g_{n} \circ g_{m} = g_{nm}$, these
functions commute. An extraordinary fact, due to Ritt
\cite{ritt1}, is that these are all such commuting rational maps.
The functions $R_{n}$ are a special class of the $g_{n}$. See
\cite{mckmoll}, section $2.13$, for details.
\end{Note}

The identity (\ref{relation1}) permits the
explicit evaluation of the zeros of $P_{m}$ and $Q_{m}$.

\begin{Prop}
\label{rootsqm}
The polynomials $P_{m}$ and $Q_{m}$ have simple real
zeros. Those of $P_{m}$ are
given by
\ba
p_{k} & = & \cot \left( \frac{(2k+1)\pi}{2m} \right)
\quad \text{ for } 0 \leq k \leq m-1, \nn
\ea
\no
and those of $Q_{m}$ are
\ba
q_{k} & = & \cot \left( \frac{k \pi}{m} \right)
\quad \text{ for } 1 \leq k \leq m-1. \nn
\ea
\end{Prop}
\begin{proof}
The identity $R_{m} = M^{-1} \circ f_{m} \circ M$ yields
\ba
R_{m}(q_{k}) & = & M^{-1}f_{m}(M( \cot(k \pi/m))) \nn \\
             & = & M^{-1}( f_{m}( e^{2 k \pi i /m})) \nn \\
             & = & M^{-1}(1) = \infty, \nn
\ea
\no
so that $Q_{m}(q_{k})=0$. The degree of $Q_{m}$ is $m-1$ and the $q_{k}$ are
all distinct, hence these are all the zeros. The argument for $p_{k}$ is
similar.
\end{proof}

\medskip

The polynomials
\ba
P^{*}_{m}(a) & = & \sum_{i=0}^{\lfloor{m/2 \rfloor}} (-1)^{i} \binom{m}{2i}
x^{2i} ~~\mbox{and}\label{pmstar} \\
Q^{*}_{m}(a) & = &
\sum_{i=0}^{\lfloor{(m-1)/2 \rfloor}} (-1)^{i} \binom{m}{2i+1}
x^{2i+1} \label{qmstar}
\ea
\no
have appeared in our development of definite integrals related to
the Hurwitz zeta function. See \cite{boesmo} for details. They are connected
to $P_{m}$ and $Q_{m}$ via
\ba
P_{m}^{*}(x) = x^{m}P_{m}(x^{-1}) & \text{ and } &
Q_{m}^{*}(x) = x^{m-1}Q_{m}(x^{-1}).
\ea

\no
The role of these polynomials in the development of the Landen
transformation comes from their trigonometric properties.

\begin{Prop}
The polynomials $P^{*}_{m}$ and $Q^{*}_{m}$ satisfy
\ba
P^{*}_{m}(\tan \theta) = \frac{\cos m \theta }{\cos^{m} \theta} & \text{ and }
& Q^{*}_{m}(\tan \theta) = \frac{\sin m \theta }{\cos^{m} \theta}.
\label{eqnpm}
\ea
\end{Prop}
\begin{proof}
We give the details for $Q_{m}^{*}$. The series expansion
 \ba
\frac{\sin(x \, \tan^{-1} t)}{( 1 + t^{2})^{x/2}} & = &
\sum_{k=0}^{\infty} \frac{(-1)^{k} \, (x)_{2k+1} }{(2k+1)!}
t^{2k+1}, \label{series1}
\ea
\no
where $(x)_{k}  =  x(x+1)(x+2)
\cdots (x+k-1)$ is the Pochhammer symbol, is established by
checking that both sides satisfy the equation \ba ( 1 + t^{2})
\frac{d^{2}g}{dt^{2}} + 2t(x+1) \frac{dg}{dt} + x(x+1) g & = & 0,
\ea \no with the initial conditions $g(0) = 0, \, g'(0) = x$. Then
$(x)_{k}$ reduces to \ba (-m)_{n} & = & (-1)^{n} n! \binom{m}{n}
\ea \no for $n \leq m$, and vanishes for $n > m$, since $m$ is an
integer. Thus (\ref{series1}) reduces to \ba \sin( m \,
\tan^{-1}t) & = & (1+t^{2})^{-m/2} \sum_{k=0}^{\lfloor{(m-1)/2
\rfloor}} (-1)^{k} \binom{m}{2k+1} t^{2k+1} \nn \ea \no for $x =
-m$. This is equivalent to the second formula in (\ref{eqnpm}). A
similar argument establishes the expression for $P_{m}^{*}$.
\end{proof}

In terms of the original polynomials, (\ref{eqnpm}) becomes
\ba
P_{m}(\cot \theta) = \frac{\cos m \theta}{\sin^{m}\theta} & \text{ and } &
Q_{m}(\cot \theta) = \frac{\sin m \theta}{\sin^{m}\theta}.
\label{rules}
\ea

\medskip

\section{A trigonometric reduction} \label{sec-trigred1}
\setcounter{equation}{0}

The example described in Section \ref{sec-example} can be extended by using
the transformation $y = R_{m}(x)$ with higher values of $m$. The explicit
evaluation of the new integrals requires knowledge of the branches of the
inverse map $x = R_{m}^{-1}(y)$. This is impractical for $m \geq 3$.
An alternative method is described in the next section.

The explicit formula for the Landen transformation
uses an expression of $ \sin^{a} \theta \, \cos^{b} \theta$, for $a, \, b
\in \mathbb{N}$, as a linear combination of trigonometric functions of
multiple angles.

We introduce the notation
\begin{equation}
c  =  \lceil{ \tfrac{a+b}{2} \rceil} \text{ and } 
d  =  \lfloor{ \tfrac{a}{2} \rfloor}.
\end{equation}
\no
The reduction formulas given below are expressed in terms of the function
\ba
T_{x}(a,b) & = & \sum_{j=0}^{x} (-1)^{a - x +j}
\binom{a}{x - j} \binom{b}{j}.
\ea
\no
Some of the identities presented here can be found in the table
appearing in \cite{gr}, page 30.

\begin{Prop}
\label{prodsincos} Let $a, \, b \in \mathbb{N}$ and $u \in
\mathbb{R}$. Then $\sin^{a}u \, \cos^{b}u$ is given by
\begin{align} \label{T-values}
& \frac{(-1)^{d}}{2^{a+b}} \, \left[ T_{c}(a,b) +
\sum_{j=1}^{c} \left( T_{c + j}(a,b)  + T_{c-j}(a,b) \right)
\cos( 2ju)  \right] &\text{ for } a
\text{ even and } b \text{ even}, \nn \\
& \frac{(-1)^{d}}{2^{a+b}}  \,
\left[ \sum_{j=1}^{c} \left( T_{c-1+j}(a,b) + T_{c-j}(a,b) \right)
\cos( (2j-1)u)  \right]  &\text{ for } a \text{ even and }
b \text{ odd}, \nn \\
& \frac{(-1)^{d}}{2^{a+b}}  \,
\left[ \sum_{j=1}^{c} \left( T_{c-1+j}(a,b) - T_{c-j}(a,b) \right)
\sin( (2j-1)u)  \right]  &\text{ for } a \text{ odd and }
b \text{ even}, \nn \\
& \frac{(-1)^{d}}{2^{a+b}}  \,
\left[ \sum_{j=1}^{c} \left( T_{c+j}(a,b) - T_{c-j}(a,b) \right)
\sin(2ju)  \right]  &\text{ for } a \text{ odd and }
b \text{ odd}. \nn
    \notag
\end{align}
\end{Prop}
\begin{proof}
Start with
\ba
(e^{iu} - e^{-iu})^{a} \, (e^{iu}+e^{-iu})^{b} & = &
\left( \sum_{k=0}^{a} \binom{a}{k} (-1)^{a-k} e^{iu(2k-a)} \right)
\left( \sum_{j=0}^{b} \binom{b}{k}  e^{iu(2j-b)} \right) \nn \\
& = & \sum_{k=0}^{a} \sum_{j=0}^{b} (-1)^{a-k} \binom{a}{k} \binom{b}{j}
e^{iu [ 2(k+j) - (a+b) ]}. \nn
\ea
\no
Therefore
\ba
\sin^{a}u \, \cos^{b}u & = & \frac{i^{-a}}{2^{a+b}}
\sum_{k=0}^{a+b} \sum_{\nu=0}^{a+b} (-1)^{a- \nu+k}
\binom{a}{\nu-k} \binom{b}{k}
e^{iu(2 \nu -a-b)}. \label{sum3}
\ea
\no
The result follows now by eliminating the imaginary terms on the right
hand side of (\ref{sum3}).
\end{proof}


\medskip

\section{The scaling of the integrand} \label{sec-scaling}
\setcounter{equation}{0}

In this section we describe a construction of the polynomials
$Z(x)$ and $E(x)$, introduced in Section \ref{sec-example}. These
are used to produce an appropriate scaling of the integrand in \ba
I & = & \ifft \frac{B(x)}{A(x)} \, dx, \label{integral1} \ea \no
so that the new denominator is $E(x)$. Recall that $E$ is the
homogeneous polynomial in the variables $(P_{m}(x), Q_{m}(x))$
introduced in Section \ref{sec-poly1}.

We now express the coefficients of $E$ and $Z$ in terms of those of $A$.
This requires the elementary symmetric functions
$$\sigma_{l}^{(p)} = \sigma_{l}^{(p)}(y_{1}, \cdots, y_{p})$$
\no
of the $p$ variables
$y_{1}, \cdots, y_{p}$. These are defined by the identity
\ba
\prod_{l=1}^{p} (y - y_{l}) & = & \sum_{l=0}^{p} (-1)^{l}
\sigma_{l}^{(p)}(y_{1}, \cdots,
y_{p}) y^{p-l}.
\label{symm-def}
\ea

\medskip

\begin{Thm}
\label{thm-scaling} Let $p, m \in \mathbb{N}$, and \ba A(x) & = &
\sum_{k=0}^{p} a_{k}x^{p-k} \ea \no be a polynomial with real
coefficients. Then there exist $mp+1$ coefficients, \ba z_{0}, \,
z_{1}, \cdots, z_{r}; \, e_{1}, \, e_{2}, \cdots, e_{p}, \ea \no
with $r = p(m-1)$, such that \ba A(x) Z(x) & = & E(x),
\label{linear} \ea \no where \ba Z(x) := \sum_{k=0}^{r}
z_{k}x^{r-k} & \text{ and } & E(x) := \sum_{l=0}^{p} e_{l} \left[
P_{m}(x) \right]^{p-l} \left[ Q_{m}(x) \right]^{l}. \ea \no The
coefficients $e_{l}$ are polynomials in the coefficients $\,
\frac{a_{1}}{a_{0}}, \cdots, \frac{a_{p}}{a_{0}}$.
\end{Thm}

\begin{Note}
The effect of the theorem is to scale the integrand $B(x)/A(x)$ to
$C(x)/E(x)$, where $C(x) = B(x)Z(x)$ and $E(x) = A(x)Z(x)$. The
degrees are recorded here: \ba \text{deg}(A) = p, \, \text{deg}(B)
= p-2, \,
\text{deg}(Z) = r = pm-p, \label{degrees} \\
\text{deg}(C) = s = pm-2, \,
\text{and deg}(E) = pm. \nn
\ea
\end{Note}

\medskip

\begin{proof}
Let $\{ x_{1}, \, x_{2}, \, \ldots, x_{p} \}$ be the roots of $A$, each
written according to its multiplicity, so that
\begin{equation}
A(x) = a_{0} \prod_{j=1}^{p} (x - x_{j}).
\end{equation}
The rational function
\begin{equation}
R_{m}(x) = \frac{P_{m}(x)}{Q_{m}(x)},
\end{equation}
\no
introduced in Proposition \ref{rofm}, is well-defined at all the roots $x_{j}$.
This follows from the fact that the roots of $Q_{m}$ are real and
our assumption that the roots
$x_{1}, \cdots, x_{p}$ of $A(x) = 0$ are not.
For $0 \leq l \leq p$, define
\begin{equation}
e_{l} := a_{0}^{m} (-1)^{l}
\prod_{j=1}^{p} Q_{m}(x_{j}) \times
\sigma_{l}^{(p)} \left( R_{m}(x_{1}), \, R_{m}(x_{2}), \ldots, R_{m}(x_{p})
\right)
\label{def-el}
\end{equation}
\no
and the polynomial
\begin{equation}
H(x) = \sum_{l=0}^{p} e_{l}x^{p-l}.
\label{def-H}
\end{equation}

We now consider the identity,
\begin{equation}
\prod_{j=1}^{p} ( y - R_{m}(x_{j}) ) =
\sum_{l=0}^{p} (-1)^{l} \sigma_{l}^{(p)} \left( R_{m}(x_{1}),
\ldots, R_{m}(x_{p}) \right) \, y^{p-l},
\end{equation}
\no
that comes from (\ref{symm-def}). Clearing  denominators, we obtain
\begin{equation}
\prod_{j=1}^{p} \left( Q_{m}(x_{j})y - P_{m}(x_{j}) \right) =
a_{0}^{-m} \sum_{l=0}^{p} e_{l}y^{p-l} = a_{0}^{-m} H(y).
\label{one}
\end{equation}
\no
In particular,
\begin{equation}
H(R_{m}(x)) = a_{0}^{m} \prod_{j=1}^{p} Q_{m}(x_{j}) \times
\prod_{j=1}^{p} \left( R_{m}(x) - R_{m}(x_{j}) \right).
\label{two}
\end{equation}
\no
Finally, define the polynomial
\ba
E(x) & = & \sum_{l=0}^{p} e_{l} P_{m}^{p-l}(x) Q_{m}^{l}(x) \nn \\
   & = & H(R_{m}(x)) Q_{m}^{p}(x). \nn
\ea
\no
Identity (\ref{two}) shows that the zeros of $E$ are precisely the
values $R_{m}(x_{j}), \, 1 \leq j \leq p$. The coefficients of $E$, given in
(\ref{one}), are symmetric polynomials of the roots $x_{j}$ of $A$. The
fundamental theorem of symmetric polynomials \cite{artin2} states that
$e_{l}$ is a polynomial in $\frac{a_{1}}{a_{0}}, \cdots, \frac{a_{p}}{a_{0}}$.
This, in turn,  proves that $e_{l} \in \mathbb{R}$ and
thus $E \in \mathbb{R}[x]$.

Now observe that (\ref{two}) yields $E(x_{j}) = 0$ and the corresponding
factor $R_{m}(x) - R_{m}(x_{j})$ appears with the same multiplicity as $x_{j}$.
We conclude that $A$ divides $E$ and define $Z$ to be the quotient. \\
\end{proof}

\medskip

\section{The reduction of the integrand} \label{sec-reduction}
\setcounter{equation}{0}

In this section, we produce explicit formulas for rational Landen
transformations of the integral \ba I & = & \ifft
\frac{B(x)}{A(x)} \, dx. \label{orig-int} \ea \no As before, we
assume that $A, \, B \in \mathbb{R}[x]$, relatively prime, and
that $I < \infty$.

In Section \ref{sec-scaling} we have scaled the integrand in (\ref{orig-int})
to the form
\ba
I & = & \ifft \frac{C(x)}{E(x)} \, dx,
\ea
\no
where the denominator is written as
\ba
E(x) & = & \sum_{l=0}^{p} e_{l} P_{m}^{p-l}(x) Q_{m}^{l}(x).
\ea
\no
Here, $P_{m}$ and $Q_{m}$ are the polynomials discussed in Section
\ref{sec-poly1}. The scaling of the denominator is achieved
through multiplication
by $Z(x)$, as  given in Theorem \ref{thm-scaling}.

The numerator becomes \ba C(x) & = & B(x)Z(x)  = \sum_{k=0}^{s}
c_{k}x^{s-k}. \ea \no The coefficients $c_{k}$ are given by \ba
c_{j} & = & \sum_{k=0}^{j} z_{k}b_{j-k} \quad \text{ for } 0 \leq
j \leq s, \ea \no where $b_{i} = 0$ if $i > p-2$ and $z_{i} =0$ if
$i > r=pm-p$.

\medskip

The change of variables $x  = \cot \theta$ and the relations (\ref{rules})
yield
\ba
I & = & \int_{0}^{\pi} \frac{\text{CT}_{m,p}(\theta) }{\text{ET}_{m,p}(\theta)}
d \theta,
\label{befsplit}
\ea
\no
where
\ba
\text{CT}_{m,p}(\theta) & = & \sum_{k=0}^{s} c_{k} \cos^{s-k}\theta \, \sin^{k}\theta
\ea
\no
and
\ba
\text{ET}_{m,p}(\theta) & = & \sum_{l=0}^{p} e_{l} \cos^{p-l}(m \theta)
\, \sin^{l}(m \theta).
\ea
\no

The discussion of this integral is divided according to the parity
of $m$. Recall that $p$ is assumed to be even. The details are
presented in the case $m$ odd.

The parameter $s = mp-2$ is even and we write $ s = 2 \lambda$. Split
(\ref{befsplit}) as
\ba
I & = & \sum_{j=0}^{\lambda} c_{2j} \int_{0}^{\pi}
\frac{\sin^{2j}\theta \, \cos^{s-2j}\theta \, d \theta }
{\text{ET}_{m,p}(\theta)}
+ \sum_{j=0}^{\lambda-1} c_{2j+1} \int_{0}^{\pi}
\frac{\sin^{2j+1}\theta \, \cos^{s-2j-1}\theta \, d \theta}
{\text{ET}_{m,p}(\theta)} \nn \\
 & \equiv & I_{1} + I_{2}, \nn
\ea
\no
and consider the evaluation of each of these integrals.  \\

\no
{\em The evaluation of $I_{1}$}. The identity in Proposition \ref{prodsincos}
yields
\ba
\sin^{2j}\theta \, \cos^{s-2j}\theta & = &
\frac{(-1)^{j}}{2^{s}}  T_{\lambda}(2j,s-2j) +  \nn \\
 & + & \frac{(-1)^{j}}{2^{s}}
\sum_{k=1}^{\lambda} \left[ T_{\lambda+k}(2j,s-2j) +
T_{\lambda-k}(2j,s-2j) \right]
\cos( 2 k \theta ), \nn
\ea
\no
and replacing this in the definition of $I_{1}$  yields
\ba
& & \label{intmess1} \\
I_{1} & = & \frac{1}{2^{s}} \sum_{j=0}^{\lambda} (-1)^{j}c_{2j}
T_{\lambda}(2j,s-2j)
\int_{0}^{\pi} \frac{d \theta}{\text{ET}_{m,p}(\theta)} + \nn  \\
 & + & \frac{1}{2^{s}} \sum_{j=0}^{\lambda}(-1)^{j}c_{2j}
\sum_{k=1}^{\lambda} \left[ T_{\lambda+k}(2j,s-2j) +
T_{\lambda-k}(2j,s-2j) \right]
\int_{0}^{\pi} \frac{\cos( 2 k \theta) \, d \theta }{\text{ET}_{m,p}(\theta)}.
\nn
\ea
\no
The periodicity of the integrand, and the fact that $mp$ is even, show that the
integrals appearing in (\ref{intmess1}) are half of the
corresponding ones over the whole period $[0, \, 2 \pi]$. Thus,
\ba
& & \label{intmess2} \\
I_{1} & = & \frac{1}{2^{s+1}} \sum_{j=0}^{\lambda} (-1)^{j}c_{2j}
T_{\lambda}(2j,s-2j)
\int_{0}^{2 \pi} \frac{d \theta}{\text{ET}_{m,p}(\theta)} + \nn  \\
 & + & \frac{1}{2^{s+1}} \sum_{j=1}^{\lambda}(-1)^{j}c_{2j}
\sum_{k=1}^{\lambda} \left[ T_{\lambda+k}(2j,s-2j) +
T_{\lambda-k}(2j,s-2j) \right]
\int_{0}^{2 \pi} \frac{\cos( 2 k \theta) \, d \theta }{\text{ET}_{m,p}(\theta)}.
\nn
\ea

\medskip

We now show that most of the integrals in (\ref{intmess2}) vanish. This
reduction is responsible for the existence of a rational Landen
transformation. \\

Introduce the notation
\ba
S_{m,p}(k) & := & \int_{0}^{2 \pi} \frac{ \sin( k \theta) \, d \theta}
{\text{ET}_{m,p}( \theta )}  \label{smp}
\ea
\no
and
\ba
C_{m,p}(k) & := & \int_{0}^{2 \pi} \frac{ \cos( k \theta) \, d \theta}
{\text{ET}_{m,p}( \theta )}.  \label{cmp}
\ea

\begin{Lem}
\label{survive}
Let $k, \, m, \, p \in \mathbb{N}$ be arbitrary. Then $S_{m,p}(k)$ and
$C_{m,p}(k)$ vanish unless $k$ is a multiple of $m$.
\end{Lem}
\begin{proof}
In the definition of $S_{m,p}(k)$ let $\theta \mapsto \theta + 2
\pi j/m$ for $j = 0, \, 1, \cdots, m-1$. The average of these $m$
integrals is \ba S_{m,p}(k) & = & \frac{1}{m} \int_{0}^{2 \pi}
\frac{d \theta}{\text{ET}_{m,p}(\theta)} \left( \sin( k \theta)
\sum_{j=0}^{m-1} \cos( 2 \pi k j/m) + \cos( k \theta)
\sum_{j=0}^{m-1} \sin( 2 \pi k j/m) \right). \nn \ea \no If $k$ is
not a multiple of $m$ the integrand vanishes because the sums in
it are the real and imaginary parts of \ba \sum_{j=0}^{m-1} e^{ 2
\pi i j k/m} & = & \frac{1 - e^{2 \pi i k}}{1 - e^{2 \pi i k/m}} =
0. \nn \ea

A similar proof follows for $C_{m,p}(k)$.
\end{proof}

\medskip

In view of Lemma \ref{survive}, we replace $k$ by $\alpha m$,
where $1 \leq \alpha \leq \nu -1$, with $\nu = p/2$. Then
(\ref{intmess2}) becomes

\ba
& & \label{intmess3} \\
I_{1} & = & \frac{1}{2^{s+1}} \sum_{j=0}^{\lambda} (-1)^{j}c_{2j}
T_{\lambda}(2j,s-2j)
\int_{0}^{2 \pi} \frac{d \theta}{\text{ET}_{m,p}(\theta)} + \nn  \\
 & + & \frac{1}{2^{s+1}} \sum_{j=1}^{\lambda}(-1)^{j}c_{2j}
\sum_{\alpha=1}^{\nu -1} \left[ T_{\lambda+ \alpha m}(2j,s-2j) +
T_{\lambda-\alpha m}(2j,s-2j) \right] \int_{0}^{2 \pi} \frac{\cos(
2 \alpha m  \theta) \, d \theta } {\text{ET}_{m,p}(\theta)}. \nn
\ea The change of variables $\varphi = m \theta$ produces \ba
\int_{0}^{2 \pi} \cdots \, d \theta = \frac{1}{m} \int_{0}^{2 \pi
m} \cdots \, d \varphi = \int_{0}^{2 \pi} \cdots \, d \varphi, \ea
\no using the periodicity of the integrand. We conclude that

\ba
& & \label{intmess4} \\
I_{1} & = & \frac{1}{2^{s+1}} \sum_{j=0}^{\lambda} (-1)^{j}c_{2j}
T_{\lambda}(2j,s-2j)
\int_{0}^{2 \pi} \frac{d \theta}{\text{ET}_{1,p}(\theta)} + \nn  \\
 & + & \frac{1}{2^{s+1}} \sum_{j=1}^{\lambda}(-1)^{j}c_{2j}
\sum_{\alpha=1}^{\nu -1} \left[ T_{\lambda+ \alpha m}(2j,s-2j) +
T_{\lambda-\alpha m}(2j,s-2j) \right]
\int_{0}^{2 \pi} \frac{\cos( 2 \alpha  \theta) \, d \theta }
{\text{ET}_{1,p}(\theta)},
\nn
\ea
\no
where the denominator is
\ba
\text{ET}_{1,p}(\theta) & = & \sum_{l=0}^{p} e_{l}
\cos^{p-l}\theta \sin^{l}\theta.
\ea

\medskip

The next step is to bring back the domain of integration to $[0, \pi]$. The
symmetry of the integrand shows that the integral over $[\pi, 2 \pi]$ is the
same as that over $[0, \pi]$. We conclude that

\ba
& & \label{intmess5} \\
I_{1} & = & \frac{1}{2^{s}} \sum_{j=0}^{\lambda} (-1)^{j}c_{2j}
T_{\lambda}(2j,s-2j)
\int_{0}^{\pi} \frac{d \theta}{\text{ET}_{1,p}(\theta)} + \nn  \\
 & + & \frac{1}{2^{s}} \sum_{j=1}^{\lambda}(-1)^{j}c_{2j}
\sum_{\alpha=1}^{\nu -1} \left[ T_{\lambda+ \alpha m}(2j,s-2j) +
T_{\lambda-\alpha m}(2j,s-2j) \right] \int_{0}^{\pi} \frac{\cos( 2
\alpha  \theta) \, d \theta } {\text{ET}_{1,p}(\theta)}. \nn \ea

\medskip

The change of variables $ y = \cot \theta$ gives, recalling that
$\nu = p/2$, \ba \int_{0}^{\pi} \frac{d
\theta}{\text{ET}_{1,p}(\theta)} & = & \ifft (1+y^{2})^{\nu-1}
\frac{dy}{H(y)}~, \ea \no where the polynomial \ba H(y) & = &
\sum_{l=0}^{p} e_{l}y^{p-l} \ea \no was introduced in
(\ref{def-H}). The identity (\ref{rules}) is now used to change
variables in the second integral to obtain

\ba
\int_{0}^{\pi} \frac{ \cos( 2 \alpha \theta) \, d \theta }
{\text{ET}_{1,p}(\theta)} & = &
\ifft ( 1 + y^{2})^{\nu - \alpha -1} \, P_{2 \alpha}(y) \, \frac{dy}{H(y)}.
\ea

\medskip

The next step is to write $P_{2 \alpha}(y)$ in terms of $1+y^{2}$.

\begin{Lem}
The polynomial $P_{2 \alpha}(y)$ can be written as
\ba
P_{2 \alpha}(y) & = & \sum_{\beta=0}^{\alpha} (-1)^{\alpha - \beta}
\frac{\alpha}
{2 \alpha - \beta} 2^{2(\alpha - \beta)} \binom{2 \alpha - \beta}{\beta}
\, ( 1 + y^{2})^{\beta}.
\ea
\end{Lem}
\begin{proof}
Start with
\ba
P_{2 \alpha}(x) & = & \sum_{j=0}^{\alpha} (-1)^{j} \binom{2 \alpha}{2 j}
x^{2 \alpha - 2j} \nn \\
& = & \sum_{j=0}^{\alpha} (-1)^{j} \binom{2 \alpha}{2 j}
\left[ ( 1+ x^{2} ) - 1 \right]^{\alpha - j} \nn \\
 & = & \sum_{\beta=0}^{\alpha} (-1)^{\alpha - \beta}
\left( \sum_{j=0}^{\alpha-\beta} \binom{2 \alpha}{2j}
\binom{\alpha-j}{\beta} \right) \, ( 1 + x^{2})^{\beta}, \nn \ea
\no and the result follows from \ba \sum_{j=0}^{\alpha - \beta}
\binom{2 \alpha}{2j} \binom{\alpha - j}{\beta} & = &
\frac{\alpha}{2 \alpha - \beta} \binom{2 \alpha - \beta}{\beta}
2^{2(\alpha - \beta)}, \text{ for } \alpha \geq \beta.
\label{sum99} \ea


\no This sum arises as a corollary of Gauss's hypergeometric
evaluation \cite{gauss12}, \ba {_{2}F_{1}} \left[ a, b; c; 1
\right] & = & \frac{\Gamma(c) \, \Gamma(c-a-b)} {\Gamma(c-a) \,
\Gamma(c-b)} \label{valueat1}~, \ea \no valid for
$\realpart{(c-a-b)}>0$. In our case, $a = \tfrac{1}{2}- \alpha, \,
b = \beta - \alpha$ and $c = \tfrac{1}{2}$, so that $c-a-b = 2
\alpha - \beta > 0$. See \cite{andrews3}, page 66 for a proof of
(\ref{valueat1}).

\end{proof}

%
%
%

We return to the evaluation of $I_{1}$. The expression in
(\ref{intmess5}) becomes \ba I_{1} & = & \frac{1}{2^{s}}
\sum_{j=0}^{\lambda} (-1)^{j} c_{2j} T_{\lambda}(2j, s-2j)
\ifft ( 1 + y^{2})^{\nu-1} \, \frac{dy}{H(y)} + \nn \\
 & + & \frac{1}{2^{s}} \sum_{j=0}^{\lambda} (-1)^{j} c_{2j}
\sum_{\alpha=1}^{\nu-1} \left( T_{\lambda + \alpha m}(2j-s -2j) +
T_{\lambda - \alpha m}(2j-s -2j) \right) \times \nn \\
 & \times &
\sum_{\beta=0}^{\alpha} (-1)^{\alpha - \beta} 2^{2(\alpha - \beta)}
\frac{\alpha}{2 \alpha - \beta} \binom{2 \alpha - \beta}{\beta}
\ifft (1 + y^{2})^{\nu - \alpha -1 + \beta} \frac{dy}{H(y)}. \nn
\ea
\medskip

\no Expanding the powers of $1+y^{2}$, and reversing the order of
summation, leads to

\ba
I_{1} & = & \frac{1}{2^{s}} \sum_{\gamma = 0}^{\nu-1}
\binom{\nu-1}{\gamma} \sum_{j=0}^{\lambda} (-1)^{j} c_{2j} T_{\lambda}(2j,s-2j)
\ifft y^{2 \gamma} \frac{dy}{H(y)} + \nn \\
&  + & \frac{1}{2^{s}} \sum_{\gamma=0}^{\nu -2 } \left( \sum_{j=0}^{\lambda}
\sum_{\alpha = 1}^{\nu-1-\gamma} \sum_{\beta=0}^{\alpha}
M_{1}(j,\alpha,\beta;m,p)
\right) \ifft y^{2 \gamma} \frac{dy}{H(y)} +  \nn \\
&  + & \frac{1}{2^{s}} \sum_{\gamma=1}^{\nu -1 } \left( \sum_{j=0}^{\lambda}
\sum_{\alpha = \nu - \gamma}^{\nu-1}
\sum_{\beta=\alpha- \nu + \gamma+1}^{\alpha} M_{1}(j,\alpha,\beta;m,p)
\right) \ifft y^{2 \gamma} \frac{dy}{H(y)},
\nn
\ea
\no
where
\ba
M_{1}(j,\alpha,\beta;m,p) & = & (-1)^{j+\alpha-\beta} c_{2j}
\frac{2^{2(\alpha - \beta)} \, \alpha}
{2 \alpha - \beta} \, \binom{2 \alpha - \beta}{\beta} \,
\binom{\nu - \alpha -1 + \beta}{\gamma}   \nn \\
& \times &  \left[ T_{\lambda + \alpha m}(2j,s-2j) +
T_{\lambda - \alpha m}(2j,s-2j) \right].  \nn
\ea

\medskip

\no
{\em The evaluation of $I_{2}$}. A similar calculation leads to

\ba
I_{2} & = &
\frac{1}{2^{s}} \sum_{\gamma=0}^{\nu -2 } \left( \sum_{j=0}^{\lambda -1}
\sum_{\alpha = 1}^{\nu-1-\gamma} \sum_{\beta=0}^{\alpha-1}
M_{2}(j,\alpha,\beta;m,p)
\right) \ifft y^{2 \gamma+1} \frac{dy}{H(y)} +  \nn \\
&  + & \frac{1}{2^{s}} \sum_{\gamma=1}^{\nu -2 } \left( \sum_{j=0}^{\lambda-1}
\sum_{\alpha = \nu - \gamma}^{\nu-1}
\sum_{\beta=0}^{\alpha-1} M_{2}(j,\alpha,\beta;m,p)
\right) \ifft y^{2 \gamma+1} \frac{dy}{H(y)},
\nn
\ea
\no
where
\ba
M_{2}(j,\alpha,\beta;m,p) & = & (-1)^{j+ \beta} c_{2j+1}
2^{2\beta+1}
\binom{\alpha+\beta}{2 \beta+1} \binom{\nu - 2 - \beta}{\gamma}  \nn \\
& \times & \left[ T_{\lambda + \alpha m}(2j+1,s-2j-1) -
T_{\lambda - \alpha m}(2j+1,s-2j-1) \right]. \nn
\ea

For the convenience of the reader we summarize the information as a theorem. \\

\bigskip

\begin{Thm}
\label{moddpeven}
Let $p, \, m \in \mathbb{N}$ and assume $p$ is even and $m$ is odd. Define
\ba
s = mp-2, \, r = p(m-1), \, \, \lambda = \frac{s}{2}, \text{ and }
\nu = \frac{p}{2},
\ea
\no
and consider the polynomials
\ba
A(x) = \sum_{k=0}^{p} a_{k}x^{p-k} & \text{ and } &
B(x) = \sum_{k=0}^{p-2} b_{k}x^{p-2-k}.
\ea
\no
Then
\ba
\ifft \frac{B(x)}{A(x)} \, dx & = &
\ifft \frac{J(x)}{H(x)} \, dx.
\ea
\no
The new denominator $H$ is given by
\ba
H(x) & = & \sum_{l=0}^{p} e_{l}x^{p-l}, \label{new-deno}
\ea
\no
where, for $0 \leq j \leq p$, the coefficients  $e_{j}$ are
solutions to the system (\ref{linear}).
An expression for $e_{j}$ in terms of the coefficients $a_{k}$ is given in
(\ref{def-el}).

The new numerator $J$ employs the function
\ba
T_{x}(a,b) & = & \sum_{j=0}^{x} (-1)^{a-x+j} \binom{a}{x-j} \binom{b}{j},
\ea
\no
and is given by
\ba
J(x) & = & \frac{1}{2^{s}} \sum_{\gamma=0}^{\nu-1}
\left(\binom{\nu-1}{\gamma} \sum_{j=0}^{\lambda}
(-1)^{j} c_{2j} T_{\lambda}(2j,s-2j)\right) x^{2 \gamma}  \label{new-nume} \\
& + & \frac{1}{2^{s}} \sum_{\gamma=0}^{\nu -2 } \left( \sum_{j=0}^{\lambda}
\sum_{\alpha = 1}^{\nu-1-\gamma} \sum_{\beta=0}^{\alpha}
M_{1}(j,\alpha,\beta;m,p) \right) x^{2 \gamma}  \nn \\
&  + & \frac{1}{2^{s}} \sum_{\gamma=1}^{\nu -1 } \left( \sum_{j=0}^{\lambda}
\sum_{\alpha = \nu - \gamma}^{\nu-1}
\sum_{\beta=\alpha- \nu + \gamma+1}^{\alpha} M_{1}(j,\alpha,\beta;m,p)
\right) x^{2 \gamma}   \nn \\
& + & \frac{1}{2^{s}} \sum_{\gamma=0}^{\nu -2 } \left( \sum_{j=0}^{\lambda -1}
\sum_{\alpha = 1}^{\nu-1-\gamma} \sum_{\beta=0}^{\alpha-1}
M_{2}(j,\alpha,\beta;m,p) \right) x^{2 \gamma +1} \nn \\
&  + & \frac{1}{2^{s}} \sum_{\gamma=1}^{\nu -2 } \left( \sum_{j=0}^{\lambda-1}
\sum_{\alpha = \nu - \gamma}^{\nu-1}
\sum_{\beta=0}^{\alpha-1} M_{2}(j,\alpha,\beta;m,p)
\right) x^{2 \gamma+1}.  \nn
\ea
\no
The coefficients $c_{j}$ are given by
\ba
c_{j} & = & \sum_{k=0}^{j} z_{k}b_{j-k} \quad \text{ for } 0 \leq j \leq s,
\ea
\no
with $b_{i} = 0$ if $i > p-2$ and $z_{i} = 0$ if $i > r = mp-p$. The
values of $z_{j}$ are obtained as solutions of the system
(\ref{linear}). \\

Finally,
\ba
M_{1}(j,\alpha,\beta;m,p) & = & (-1)^{j+\alpha-\beta} c_{2j}
\frac{2^{2(\alpha - \beta)} \, \alpha}
{2 \alpha - \beta} \, \binom{2 \alpha - \beta}{\beta} \,
\binom{\nu - \alpha -1 + \beta}{\gamma}   \nn \\
& \times &  \left[ T_{\lambda + \alpha m}(2j,s-2j) +
T_{\lambda - \alpha m}(2j,s-2j) \right],  \nn
\ea
\no
and
\ba
M_{2}(j,\alpha,\beta;m,p) & = & (-1)^{j+ \beta} c_{2j+1}
2^{2\beta+1}
\binom{\alpha+\beta}{2 \beta+1} \binom{\nu - 2 - \beta}{\gamma}  \nn \\
& \times & \left[ T_{\lambda + \alpha m}(2j+1,s-2j-1) -
T_{\lambda - \alpha m}(2j+1,s-2j-1) \right]. \nn
\ea
\end{Thm}

\bigskip

\begin{Note}
Surprisingly, the expressions for $H$ and $J$ given in
(\ref{new-deno}) and (\ref{new-nume}) remain valid when $m$ is
even.  This results from a similar calculation whose details are
omitted here.
\end{Note}

\bigskip

\section{Examples of rational Landen transformations} \label{sec-landenex}
\setcounter{equation}{0}

This section contains some examples that illustrate the rational Landen
transformations.  \\

\begin{Example}
We calculate the transformation
for the case $p = m = 2$. The integrand in 
\ba I & = & \ifft
\frac{dx}{a_{0}x^{2} + a_{1}x + a_{2}}, \text{ with } a_{0} \neq
0~, \label{example2} \ea 
\no
is quadratic, thus $p=2$.
We construct the Landen transformation with convergence order $m =
2$. Therefore, $s = mp-2 = 2$ and $r = p(m-1) = 2$ in the notation
defined in Theorem \ref{thm-scaling}. The scaling of Section
\ref{sec-scaling} amounts to finding parameters $z_{0}, z_{1},
z_{2},$ and $e_{0}, e_{1},e_{2}$ such that 
\ba 
(a_{0}x^{2} + a_{1}x
+ a_{2})(z_{0}x^{2}+z_{1}x+z_{2}) & = & e_{0}P_{2}^{2}(x) +
e_{1}P_{2}(x)Q_{2}(x) + e_{2}Q_{2}^{2}(x), \nn 
\ea 
\no 
with $P_{2}(x) = x^{2}-1$ and $Q_{2}(x) = 2x$. The linear system
(\ref{linear}) is of order $mp+1 = 5$ and we choose the free
parameter 
\ba 
e_{0} & = & Q_{2}(x_{1}) Q_{2}(x_{2})  = 4x_{1}x_{2}
= \frac{4a_{2}}{a_{0}}, 
\ea 
\no 
according to the convention in
(\ref{def-el}). The solution of (\ref{linear}) yields \ba z_{0} =
\frac{4 a_{2}}{a_{0}^{2}}, \,\,\,  z_{1} =
-\frac{4a_{1}}{a_{0}^{2}}, \,\, \, z_{2}= \frac{4}{a_{0}}, \ea \no
and the formulas in (\ref{def-el}) produce 
\ba 
e_{0}=
\frac{4a_{2}}{a_{0}}, \,\,\, e_{1} = \frac{2a_{1}(a_{2}-a_{0})}
{a_{0}^{2}}, \,\,\,  e_{2} =
\frac{(a_{0}+a_{2})^{2}-a_{1}^{2}}{a_{0}^{2}}. 
\ea 
\no 
Therefore,
the denominator of the new integrand is 
\ba 
H(x) & = &
\frac{4a_{2}}{a_{0}}x^{2} + \frac{2a_{1}(a_{2}-a_{0})}{a_{0}^{2}}x
+ \frac{(a_{0}+a_{2})^{2} - a_{1}^{2}}{a_{0}^{2}}. 
\ea
 \no 
The new numerator is obtained from the formulas given in Theorem
\ref{moddpeven}. In this case $\lambda = 1$ and $\nu = 1$, thus
only one sum contributes to its value: 
\ba 
J(x) & = & \frac{2(a_{0}+a_{2})}{a_{0}^{2}}. 
\ea
 \no 
We conclude that 
\ba
\ifft \frac{dx}{a_{0}x^{2}+a_{1}x + a_{2}} & = & \ifft
\frac{2(a_{0}+a_{2}) \, dx}{ 4a_{0}a_{2}x^{2} +
2a_{1}(a_{2}-a_{0})x + \left[ (a_{0}+a_{2})^{2} - a_{1}^{2}
\right]}, \nn 
\ea
\no
and (\ref{example2}) is 
invariant under the transformation 
\ba
a_{0} & \mapsto & \frac{2a_{0}a_{2}}{a_{0}+a_{2}}, \label{map1} \\
a_{1} & \mapsto & \frac{a_{1}(a_{2}-a_{0})}{a_{0}+a_{2}}, \nn \\
a_{2} & \mapsto & \frac{(a_{0}+a_{2})^{2} -
a_{1}^{2}}{2(a_{0}+a_{2})}. \nn 
\ea
\no
This was announced in (\ref{map0}).
\end{Example}

\bigskip

\begin{Example}
We present the Landen transformation of order $3$ for the rational
function \ba R(x) & = & \frac{x^{2}+4x+4} {x^{6} + 16x^{5} +
114x^{4} + 452x^{3}+1041x^{2}+1300x+676}. \ea \no This is an
example that violates the main assumption on the nature of the
roots of $A$. Indeed, \ba A(x) = (x+2)^{2} (x^{2}+6x+13)^{2}  &
\text{ and } & B(x) = (x+2)^{2}, \ea \no so that $A$ has real
roots. Although $R$ is not reduced, it is integrable
over $\mathbb{R} \, :$ 
\ba I & = & \ifft
\frac{dx}{(x^{2}+6x+13)^{2}} = \frac{\pi}{16}. \nn \ea \no This
example shows that the rational Landen transformations preserve
the existence of real poles of the integrand. Moreover, the real
zeros that cancel these singularities are transformed accordingly
to preserve convergence.

The roots of $A(x)=0$ are \ba x_{1} = x_{2} = -3-2i, \,
x_{3}=x_{4}=-3+2i, \, x_{5}=x_{6} = -2. \ea \no The value of
$e_{0}$ given in (\ref{def-el}) yields \ba e_{0} & = &
\prod_{k=1}^{6} Q_{3}(x_{k}) = 269353744. \ea \no We have that
$p=6$ and $m=3$, and so $r=12$ and $s=16$. Solving a system of
order $19$ yields \ba A_{1}(x) & = & (11x+2)^{2} ( 373x^{2} + 594x
+ 481)^{2} \ea \no as the new denominator, and  the new numerator
is \ba B_{1}(x) & = & (11x+2)^{2}( 854x^{2} + 3240x + 10709). \ea
\no Observe that the algorithm preserves the existence of a real
root, but the root at $x = -2/11$ is cancelled. The reader will
check the invariance:  \ba \ifft \frac{854x^{2} + 3240x +
10709}{(373x^{2} + 594x + 481)^{2}} \, dx & = & \frac{\pi}{16}~.
\ea \no
\end{Example}

\bigskip

\begin{Example}
Finally we present a numerical example that illustrates the convergence of
the iterative transformations constructed in this paper.  The
original integral is written in the form \ba I & = &
\frac{b_{0}}{a_{0}} \ifft \frac{x^{p-2} + b_{0}^{-1}b_{1}x^{p-3} +
b_{0}^{-1}b_{2}x^{p-4} + \cdots + b_{0}^{-1}b_{p-2}} {x^{p} +
a_{0}^{-1}a_{1}x^{p-1} + a_{0}^{-1}a_{2}x^{p-2} + \cdots +
a_{0}^{-1}a_{p}} \, dx. \ea \no The Landen transformation
generates a sequence of coefficients, \ba {\mathfrak{P}}_{n} & :=
& \{ a_{0}^{(n)}, \, a_{1}^{(n)}, \cdots, a_{p}^{(n)}; \,
b_{0}^{(n)}, \, b_{1}^{(n)}, \cdots, b_{p-2}^{(n)} \, \}~, \ea \no
with ${\mathfrak{P}}_{0} = \mathfrak{P}$ as in
(\ref{initialpara}). We wish to show that, as $n \to \infty$,
\begin{equation}
u_{n} := \left( \frac{a_{1}^{(n)}}{a_{0}^{(n)}},
\frac{a_{2}^{(n)}}{a_{0}^{(n)}},  \cdots,
\frac{a_{p}^{(n)}}{a_{0}^{(n)}},
\frac{b_{1}^{(n)}}{b_{0}^{(n)}},
\frac{b_{2}^{(n)}}{b_{0}^{(n)}},  \cdots,
\frac{b_{p-2}^{(n)}}{b_{0}^{(n)}} \right)
\end{equation}
\no
converges to
\begin{equation}
u_{\infty} := \left( 0, \binom{q}{1},0, \binom{q}{2}, \cdots,
\binom{q}{q}; 0, \binom{q-1}{1},0, \binom{q-1}{2}, \cdots,
\binom{q-1}{q-1} \right)~,
\end{equation}
\no
where $q = p/2$. The
invariance of the integral then shows that
\ba
\frac{b_{0}^{(n)}}{a_{0}^{(n)}} \to \frac{1}{\pi} I.
\label{conv-int}
\ea

\medskip

The convergence of $v := u_{n}- u_{\infty}$ to $0$ is measured in
the $L_{2}-$norm, \ba \Vert v \Vert_{2} & = & \frac{1}{\sqrt{2p-2}}
\left( \sum_{k=1}^{2p-2} \| v_{k} \|^{2} \right)^{1/2}, \ea \no and
also the $L_{\infty}$-norm, \ba \Vert v \Vert_{\infty} & = &
\text{Max} \left\{ \| v_{k} \|: 1 \leq k \leq 2p-2 \, \right\}.
\ea \no The rational functions appearing as integrands have
rational coefficients, so, as a measure of their complexity, we
take the largest number of digits of these coefficients. This
appears in the column marked {\em size}.

\bigskip

The following tables illustrate the iterates of rational Landen
transformations of order $2, \, 3$ and $4$, applied to the example
\ba I & = & \ifft \frac{3x+5}{x^{4}+ 14x^3+74x^2+184x+208} \, dx =
- \frac{7 \pi}{12}. \ea

\medskip

\no The first column gives the $L_{2}$-norm of $u_{n} -
u_{\infty}$, the second its $L_{\infty}$-norm, the third presents
the relative error in (\ref{conv-int}), and in the last column we
give the size of the rational integrand. At each step, we verify
that the new rational function integrates to $-7 \pi/12$.

\medskip

\begin{center}
Method of order $2$
\end{center}

\medskip

\begin{center}
\begin{tabular}{||c|c|c|c|c||}
\hline
$n$ & $L_{2}$-norm & $L_{\infty}$-norm & Error & Size \\ \hline
$1$ & $58.7171$ & $69.1000$ & $1.02060$ & $5$ \\
$2$ & $7.444927$ & $9.64324$ & $1.04473$ & $10$ \\
$3$ & $4.04691$ & $5.36256$ & $0.945481$ & $18$ \\
$4$ & $1.81592$ & $2.41858$ & $1.15092$ & $41$ \\
$5$ & $0.360422$ & $0.411437$ & $0.262511$ & $82$ \\
$6$ & $0.0298892$ & $0.0249128$ & $0.0189903$ & $164$ \\
$7$ & $0.000256824$ & $0.000299728$ & $0.0000362352$ & $327$ \\
$8$ & $1.92454 \times 10^{-8}$ & $2.24568 \times 10^{-8}$ & $1.47053 \times 10^{-8}$ & $659$ \\
$9$ & $1.0823 \times 10^{-16}$ & $1.2609 \times 10^{-16}$ & $8.2207 \times 10^{-17}$
& $1318$ \\
\hline
\end{tabular}
\end{center}

\medskip

As expected we observe quadratic convergence in the $L_{2}-$norm
and also in the $L_{\infty}-$norm. The size of the integrand is
doubled at each iteration.

\newpage

\begin{center}
Method of order $3$
\end{center}

\medskip

\begin{center}
\begin{tabular}{||c|c|c|c|c||}
\hline
$n$ & $L_{2}$-norm & $L_{\infty}$-norm & Error & Size \\ \hline
$1$ & $15.2207$ & $20.2945$ & $1.03511$ & $8$ \\
$2$ & $1.97988$ & $1.83067$ & $0.859941$ & $23$ \\
$3$ & $0.41100$ & $0.338358$ & $0.197044$ & $69$ \\
$4$ & $0.00842346$ & $0.00815475$ & $0.00597363$ & $208$ \\
$5$ & $5.05016 \times 10^{-8}$ & $5.75969 \times 10^{-8}$ & $1.64059 \times 10^{-9}$ & $626$ \\
$6$ & $1.09651 \times 10^{-23}$ & $1.02510 \times 10^{-23}$ & $3.86286 \times 10^{-24}$
& $1878$ \\
$7$ & $1.12238 \times 10^{-70}$ & $1.22843 \times 10^{-70}$ & $8.59237 \times 10^{-71}$ & $5634$  \\
\hline
\end{tabular}
\end{center}

\medskip

\begin{center}
Method of order $4$
\end{center}

\medskip

\begin{center}
\begin{tabular}{||c|c|c|c|c||}
\hline
$n$ & $L_{2}$-norm & $L_{\infty}$-norm & Error & Size \\ \hline
$1$ & $7.44927$ & $9.64324$ & $1.04473$ & $10$ \\
$2$ & $1.81592$ & $2.41858$ & $1.15092$ & $41$ \\
$3$ & $0.0298892$ & $0.0249128$ & $0.0189903$ & $164$ \\
$4$ & $1.92454 \times 10^{-8}$ & $2.249128 \times 10^{-8}$ & $1.47053 \times 10^{-8}$ & $659$ \\
$5$ & $3.40769 \times 10^{-33}$ & $3.96407 \times 10^{-33}$ & $2.56817 \times 10^{-33}$ & $2637$ \\
\hline
\end{tabular}
\end{center}

\end{Example}

\section{Conclusions} \label{sec-conclusions}
\setcounter{equation}{0}

We have presented an algorithm for the evaluation of a rational integral over
$\mathbb{R}$. Numerical evidence of its convergence is presented.

\bigskip

\no
{\bf Acknowledgments}. The work of the second author was partially funded by
$\text{NSF-DMS } 0409968$. The first author was partially supported as a
graduate student by the same grant.

\bigskip

\bibliography{../../../AllRef/small2}
\bibliographystyle{plain}
\end{document}